\documentclass[12pt]{amsart}
\usepackage{amsmath,amssymb,amsbsy,amsfonts,latexsym,amsopn,amstext,cite,
                                               amsxtra,euscript,amscd,bm}
\usepackage{booktabs}
\usepackage{url}

\usepackage{mathrsfs}

\usepackage{color}
\usepackage[colorlinks,linkcolor=blue,anchorcolor=blue,citecolor=blue,backref=page]{hyperref}

\usepackage{graphics,epsfig}
\usepackage{graphicx}
\usepackage{float}
\usepackage{epstopdf}
\hypersetup{breaklinks=true}
\usepackage{verbatim}
\usepackage{graphicx}
\usepackage{caption}
\usepackage{subcaption}

\usepackage{soul} 
\usepackage{ulem}

\usepackage{float}
\usepackage{epstopdf}
\hypersetup{breaklinks=true}

\usepackage[np]{numprint}
\npdecimalsign{\ensuremath{.}}

\usepackage{bibentry}

\usepackage[english]{babel}
\usepackage{mathtools}
\usepackage{todonotes}
\usepackage{url}
\usepackage[norefs,nocites]{refcheck}

\begin{document}
\newtheorem{theorem}{Theorem}
\newtheorem*{thm*}{Theorem}
\newtheorem{lemma}[theorem]{Lemma}
\newtheorem{lem}[theorem]{Lemma}
\newtheorem{prop}[theorem]{Proposition}
\newtheorem{cor}[theorem]{Corollary}
\newtheorem{prob}[theorem]{Problem}
\newtheorem*{conj}{Conjecture}
\newtheorem{question}[theorem]{Question}
\newtheorem{defn}[theorem]{Definition}
\newtheorem{rem}[theorem] {Remark} 



\numberwithin{equation}{section}
\numberwithin{theorem}{section}
\numberwithin{table}{section}
\numberwithin{figure}{section}

\newcommand{\rad}{\operatorname{rad}}

\newcommand{\Z}{{\mathbb Z}} 
\newcommand{\Q}{{\mathbb Q}}
\newcommand{\R}{{\mathbb R}}
\newcommand{\C}{{\mathbb C}}
\newcommand{\N}{{\mathbb N}}
\newcommand{\F}{{\mathbb F}}
\newcommand{\FF}{{\mathbb F}}
\newcommand{\fq}{\mathbb{F}_q}
\newcommand{\rmk}[1]{\footnote{{\bf Comment:} #1}}

\renewcommand{\mod}{\;\operatorname{mod}}
\newcommand{\ord}{\operatorname{ord}}
\newcommand{\TT}{\mathbb{T}} 
 
\renewcommand{\i}{{\mathrm{i}}}
\renewcommand{\d}{{\mathrm{d}}}
\renewcommand{\^}{\widehat}
\newcommand{\HH}{\mathbb H}
\newcommand{\Vol}{\operatorname{vol}}
\newcommand{\area}{\operatorname{area}}
\newcommand{\dist} {\operatorname{dist}}
\newcommand{\adj} {\operatorname{adj}}

\newcommand{\tr}{\operatorname{tr}}
\newcommand{\norm}{\mathcal N} 
\newcommand{\disc}{\operatorname{disc}}

\newcommand{\intinf}{\int_{-\infty}^\infty}

\newcommand{\ave}[1]{\left\langle#1\right\rangle} 
\newcommand{\Var}{\operatorname{Var}}
\newcommand{\Prob}{\operatorname{Prob}}
\newcommand{\Cov}{\operatorname{cov}}

\newcommand{\sym}{\operatorname{Sym}}
\newcommand{\Sym}{\operatorname{Sym}}

\newcommand{\CA}{{\mathcal C}_A}
\newcommand{\cond}{\operatorname{cond}} 
\newcommand{\lcm}{\operatorname{lcm}}
\newcommand{\Kl}{\operatorname{Kl}} 
\newcommand{\leg}[2]{\left( \frac{#1}{#2} \right)}  
\newcommand{\Li}{\operatorname{Li}}

\newcommand{\sumstar}{\sideset \and^{*} \to \sum}

\newcommand{\sumf}{\sum^\flat}
 
\newcommand{\conv}{*}

\newcommand{\Ht}{\operatorname{Ht}}

\newcommand{\E}{\operatorname{\mathbb E}} 
\newcommand{\sign}{\operatorname{sign}} 
\newcommand{\meas}{\operatorname{meas}} 
\newcommand{\length}{\operatorname{length}} 

\newcommand{\GL}{\operatorname{GL}}
\newcommand{\SL}{\operatorname{SL}}
\newcommand{\Sp}{\operatorname{Sp}} 
\newcommand{\USp}{\operatorname{USp}} 

\newcommand{\re}{\operatorname{Re}}
\newcommand{\im}{\operatorname{Im}}
\newcommand{\res}{\operatorname{Res}}
  
\newcommand{\diam}{\operatorname{diam}}
 
\newcommand{\fixme}[1]{\footnote{Fixme: #1}}
 
\newcommand{\orb}{\operatorname{Orb}}
\newcommand{\supp}{\operatorname{Supp}}

\newcommand{\OP}{\operatorname{Op}} 
\newcommand{\opS}{\operatorname{S}} 
\newcommand{\opT}{\operatorname{T}} 
\newcommand{\HN}{{\cH}_{N}}
\newcommand{\OPN}{\operatorname{Op}_N} 
\newcommand{\TN}{\opT_N} 
\newcommand{\TNr}{\opT_N^{(r)}} 
\newcommand{\UNr}{U_{N,r}} 
\newcommand{\TM}{\opT_M} 
\newcommand{\Hp}{\cH_p} 
\newcommand{\OPp}{\operatorname{Op}_p} 
\newcommand{\Tp}{\operatorname{T}_p} 
 \newcommand{\diag}{\operatorname{diag}}

 \newcommand{\Fp}{\mathbb F_p}
 \newcommand{\Gal}{\operatorname{Gal}} 

 
\newcommand{\commZ}[2][]{\todo[#1,color=orange!60]{Z: #2}}
\newcommand{\commI}[2][]{\todo[#1,color=green!60]{I: #2}}
\newcommand{\commII}[2][]{\todo[#1,color=magenta!60]{I: #2}}
\newcommand{\commP}[2][]{\todo[#1,color=red!60]{P: #2}}
\newcommand{\commA}[2][]{\todo[#1,color=yellow!60]{A: #2}}
\newcommand{\commS}[2][]{\todo[#1,color=blue!60]{S: #2}}

\def\ccr#1{\textcolor{red}{#1}}
\def\cco#1{\textcolor{orange}{#1}}
\def\ccc#1{\textcolor{cyan}{#1}}

\def\Degf{e} 

\def\cont{\mathrm{cont}}

\newcommand{\bs}{\boldsymbol}
\def \a{\alpha} \def \b{\beta} \def \d{\delta} \def \e{\varepsilon} \def \g{\gamma} \def \k{\kappa} \def \l{\lambda} \def \s{\sigma} \def \t{\vartheta} \def \z{\zeta}

\newcommand{\mb}{\mathbb}

\def\sssum{\mathop{\sum\!\sum\!\sum}}
\def\ssum{\mathop{\sum\ldots \sum}}
\def\iint{\mathop{\int\ldots \int}}

\newcommand{\cc}[1]{{\color{magenta} #1}}

\newfont{\teneufm}{eufm10}
\newfont{\seveneufm}{eufm7}
\newfont{\fiveeufm}{eufm5}
%
%
\newfam\eufmfam
     \textfont\eufmfam=\teneufm
\scriptfont\eufmfam=\seveneufm
     \scriptscriptfont\eufmfam=\fiveeufm
%
%
\def\frak#1{{\fam\eufmfam\relax#1}}

\newcommand{\bflambda}{{\boldsymbol{\lambda}}}
\newcommand{\bfmu}{{\boldsymbol{\mu}}}
\newcommand{\bfxi}{{\boldsymbol{\eta}}}
\newcommand{\bfrho}{{\boldsymbol{\rho}}}

\def\eps{\varepsilon}

\def\fI{\mathfrak I}
\def\fK{\mathfrak K}
\def\fT{\mathfrak{T}}
\def\fL{\mathfrak L}
\def\fR{\mathfrak R}
\def\fQ{\mathfrak Q}

\def\fA{{\mathfrak A}}
\def\fB{{\mathfrak B}}
\def\fC{{\mathfrak C}}
\def\fM{{\mathfrak M}}
\def\fS{{\mathfrak  S}}
\def\fU{{\mathfrak U}}

\def\cPg{\cP_{\mathrm{good}}}

\def\cNg{\cN_{\mathrm{good}}}
\def\tcNg{\widetilde\cN_{\mathrm{good}}}
 
\def\sssum{\mathop{\sum\!\sum\!\sum}}
\def\ssum{\mathop{\sum\ldots \sum}}
\def\dsum{\mathop{\quad \sum \qquad \sum}}
\def\iint{\mathop{\int\ldots \int}}
 
\def\T {\mathsf {T}}
\def\Tor{\mathsf{T}_d}
\def\Tore{\widetilde{\mathrm{T}}_{d} }

\def\sM {\mathsf {M}}
\def\sL {\mathsf {L}}
\def\sK {\mathsf {K}}
\def\sP {\mathsf {P}}

\def\ss{\mathsf {s}}

\def \vk{\vec{k}}
\def \vl{\boldsymbol{\ell}}
\def \vm{\vec{m}}
\def \vn{\vec{n}}
\def \vx{\vec{x}}
\def \vf{\vec{f}}
\def \vu{\vec{u}}
\def \vv{\vec{v}}
\def \vw{\vec{w}}
\def \vz{\vec{z}}
\def\va{\vec{a}}
\def\vb{\vec{b}}

\def \balpha{\bm{\alpha}}
\def \bbeta{\bm{\beta}}
\def \bgamma{\bm{\gamma}}
\def \bdelta{\bm{\delta}}
\def \bzeta{\bm{\zeta}}
\def \blambda{\bm{\lambda}}
\def \bchi{\bm{\chi}}
\def \bphi{\bm{\varphi}}
\def \bpsi{\bm{\psi}}
\def \bxi{\bm{\xi}}
\def \bnu{\bm{\nu}}
\def \bmu{\bm{\mu}}
\def \bomega{\bm{\omega}}

\def \bell{\bm{\ell}}

\def\eqref#1{(\ref{#1})}

\def\vec#1{\mathbf{#1}}

\newcommand{\abs}[1]{\left| #1 \right|}

\def\Zq{\mathbb{Z}_q}
\def\Zqx{\mathbb{Z}_q^*}
\def\Zd{\mathbb{Z}_d}
\def\Zdx{\mathbb{Z}_d^*}
\def\Zf{\mathbb{Z}_f}
\def\Zfx{\mathbb{Z}_f^*}
\def\Zp{\mathbb{Z}_p}
\def\Zpx{\mathbb{Z}_p^*}
\def\cM{\mathcal M}
\def\cE{\mathcal E}
\def\cH{\mathcal H}

\def\le{\leqslant}
\def\leq{\leqslant}
\def\ge{\geqslant}
\def\leq{\leqslant}

\def\sfB{\mathsf {B}}
\def\sfC{\mathsf {C}}
\def\sfS{\mathsf {S}}
\def\sfI{\mathsf {I}}
\def\sfT{\mathsf {T}}
\def\L{\mathsf {L}}
\def\FF{\mathsf {F}}

\def\sE {\mathscr{E}}
\def\sS {\mathscr{S}}
\def\Gal {\mathrm{Gal}}

\def\cA{{\mathcal A}}
\def\cB{{\mathcal B}}
\def\cC{{\mathcal C}}
\def\cD{{\mathcal D}}
\def\cE{{\mathcal E}}
\def\cF{{\mathcal F}}
\def\cG{{\mathcal G}}
\def\cH{{\mathcal H}}
\def\cI{{\mathcal I}}
\def\cJ{{\mathcal J}}
\def\cK{{\mathcal K}}
\def\cL{{\mathcal L}}
\def\cM{{\mathcal M}}
\def\cN{{\mathcal N}}
\def\cO{{\mathcal O}}
\def\cP{{\mathcal P}}
\def\cQ{{\mathcal Q}}
\def\cR{{\mathcal R}}
\def\cS{{\mathcal S}}
\def\cT{{\mathcal T}}
\def\cU{{\mathcal U}}
\def\cV{{\mathcal V}}
\def\cW{{\mathcal W}}
\def\cX{{\mathcal X}}
\def\cY{{\mathcal Y}}
\def\cZ{{\mathcal Z}}
\newcommand{\rmod}[1]{\: \mbox{mod} \: #1}

\def\cg{{\mathcal g}}

\def\vy{\mathbf y}
\def\vr{\mathbf r}
\def\vx{\mathbf x}
\def\va{\mathbf a}
\def\vb{\mathbf b}
\def\vc{\mathbf c}
\def\ve{\mathbf e}
\def\vf{\mathbf f}
\def\vg{\mathbf g}
\def\vh{\mathbf h}
\def\vk{\mathbf k}
\def\vm{\mathbf m}
\def\vz{\mathbf z}
\def\vu{\mathbf u}
\def\vv{\mathbf v}

\def\e{{\mathbf{\,e}}}
\def\ep{{\mathbf{\,e}}_p}
\def\eq{{\mathbf{\,e}}_q}
\def\ek{{\mathbf{\,e}}_k}

\def\Tr{{\mathrm{Tr}}}
\def\Nm{{\mathrm{Nm}}}

\def\newr{t}

 \def\SS{{\mathbf{S}}}

 \def\0{{\mathbf{0}}}

\def\({\left(}
\def\){\right)}
\def\l|{\left|}
\def\r|{\right|}
\def\fl#1{\left\lfloor#1\right\rfloor}
\def\rf#1{\left\lceil#1\right\rceil}
\def\sumstar#1{\mathop{\sum\vphantom|^{\!\!*}\,}_{#1}}

\def\mand{\qquad \mbox{and} \qquad}

\def\tblue#1{\begin{color}{blue}{{#1}}\end{color}}

\newcommand{\fp}{\mathfrak p}

 \def \xbar{\overline x}

\newif\ifcomment

\let\varepsilon\varepsilon

\title[Higher dimensional cat maps modulo prime powers]{On the equidistribution of eigenfunctions of higher dimensional cat maps modulo prime powers}

\author[S. Bhakta]{Subham Bhakta}
\address{School of Mathematics and Statistics, University of New South Wales, Sydney, NSW 2052, Australia.} 
\email{subham.bhakta@unsw.edu.au}

\author[I. E. Shparlinski] {Igor E. Shparlinski}
\address{School of Mathematics and Statistics, University of New South Wales, Sydney NSW 2052, Australia}
\email{igor.shparlinski@unsw.edu.au}

\begin{abstract}
The so-called {\it Cat Map\/}, which is a  linear map generated by symplectic  
matrices $A \in \operatorname{Sp}(2d,\mathbb Z)$ modulo an integer $N\ge 1$  is a classical  
discrete model of quantum ergodicity. For $d=1$ it has been investigated by P.~Kurlberg and Z.~Rudnick (2001), J.~Bourgain  (2007) and  A.~Ostafe, I.~E.~Shparlinski  and J.~F.~Voloch (2023). 
The significantly more difficult case of   $d>1$  has been studied by P.~Kurlberg, A.~Ostafe,  Z.~Rudnick and I.~E.~Shparlinski  (2025). All these works apply only to almost all $N$ and are unable to provide an explicit sequence of such ``good'' values of $N$. The authors (2026) have given 
such a construction with a polynomial decay on the corresponding matrix coefficients. 
Here we use some new ideas and improve this bound. 
\end{abstract}



\maketitle

\tableofcontents 

\section{Introduction}
\subsection{Quantised linear maps and equidistribution of their eigenfunctions}
We need to recall some basic notations and definitions, in which we follow very closely 
our previous work~\cite{BhaShp}, which in turn is based on 
the exposition in~\cite{KORS}.

Let $\Z_N$ denote the residue ring modulo  an integer $N$ 
and let  $\HN = L^2\(\Z_N^d\)$ be  the   Hilbert space   with the  scalar product
$\langle \varphi_1,  \varphi_2 \rangle$ 
and norm $ \|\varphi\|$, given by
\[
\langle \varphi_1,  \varphi_2 \rangle = \frac{1}{N^d} \sum_{\vec u \in \Z_N^d}  \varphi_1(\vec u) \overline{ \varphi_2(\vec u)} \quad  \mand \|\varphi\| = \langle \varphi ,  \varphi  \rangle^{1/2} 
\]
respectively, where  $\varphi_1, \varphi_2, \varphi \in \HN$.

We then consider  the family of  unitary operators
\[
\TN(\vu):\HN\to \HN, \qquad \vu=(\vu_1,
\vu_2)\in \Z^d\times \Z^d=\Z^{2d},
\]
which are defined by the following action on $\varphi \in \HN$
\begin{equation}
\label{eq:TN}
\left(\TN(\vu)\varphi \right)(\vec w) = \e_{2N}(\vu_1\cdot \vu_2)\e_N(\vu_2\cdot \vec w) \varphi (\vec w+\vu_1),
\end{equation}
for any $\vw \in \Z_N^d$,
where 
\[
\e(z) =  \exp\(2\pi i z\),\quad \ek(z)=\e(z/k), 
\]
and we always follow the convention that integer arguments of functions on $\Z_N$ are reduced modulo $N$.
We also see that~\eqref{eq:TN} implies  
\[
\TN(\vu) \TN(\vv)  = \e_{2N}\left(\omega\left(\vu,\vv\right)\right)\TN(\vu+\vv),
\]
 where for  $\vec x=(\vec x_1,\vec x_2)$, 
$\vec y=(\vec y_1,\vec y_2)\in \R^d\times \R^d$ we define
\begin{equation}
\label{eq:omega}
\omega(\vec x,\vec y)  = \vec x_1 \cdot \vec y_2-\vec x_2\cdot \vec y_1,
\end{equation}
see also~\cite[Equation~(2.6)]{KR2001}.

For each real-valued function $f\in C^\infty(\TT^{2d})$ (an ``observable''), where $\TT = \R/Z$ is a unit torus, we define  the following operator   on $\HN$:   \begin{equation}
\label{eq:OpN}
\OP_{N}(f) = \sum_{\vu\in \Z^{2d}}\^f(\vu) \TN(\vu),
\end{equation}
where $ \^ f(\vec u)$ are the corresponding Fourier coefficients: 
\[
f(\vx)= \sum_{\vu\in \Z^{2d}} \^ f(\vu)\e(\vu \cdot \vx).
\]

We now fix a $2d \times 2d$ symplectic matrix  $A\in \Sp(2d,\Z)$ with  $A \equiv I_{2d}\pmod 2$,  where $I_{2d}$ is the $2d$-dimensional identity matrix.
The matrix $A$ induces in a  natural way  a linear map on  $\Z_N^{2d}$, which is   
also known as the {\it Cat Map\/}.
 It is also useful to notice that $A$ preserves the symplectic form~\eqref{eq:omega}, that is, $\omega(A \vec x, A \vec y) = \omega(\vec x,\vec y)$.

We further assume that $A$ has an irreducible characteristic polynomial (and thus is diagonalisable over $\C$) and there are no roots of unity amongst the eigenvalues of $A$ and their nontrivial ratios.

As in~\cite{BhaShp,KORS}, we recall that for each $N\geq 1$, there is a
 unitary operator $U_N(A)$ on $\HN$ such that for every    $f\in C^\infty(\TT^{2d})$, we have \[
U_N(A)^*\OPN(f) U_N(A) = \OPN(f\circ A) , 
\]
where $U_N(A)^*= \overline{U_N(A)}^{\,t}$,  which is known as {\it the exact Egorov property\/}.
We refer to~\cite{GuHa, KRR,KRDuke, KR2001, KR2, Rudnicksurvey} for  a detailed exposition in the case $d=1$ and~\cite{KelmerAnnals} for higher dimensions.

 Based on general results of  the Quantum Ergodicity Theorem~\cite{Schnirelman, Zelditch, CdV}, 
one may  expect that for  any normalised 
 sequence of eigenfunctions 
 $\psi_N\in \HN$ of  the operator $U_N(A)$  we have
\[
\left \langle \OPN(f)\psi_N,\psi_N \right \rangle \to \int_{\TT^{2d}} f(\vec x)d\vec x 
\] 
for all $f\in C^\infty(\TT^{2d})$, 
as $N\to \infty$ along  a reasonably ``typical" sequence $\cN$ of integers.
In this case we say that the sequence of eigenfunctions $\{\psi_N\}_{N \in \cN}$ is equidistributed, 
   see also~\cite{RudnickSarnakQUE}. 
   
Following~\cite{BhaShp,KORS,OSV}, we now introduce the following quantitative measure of the above convergence, which 
 we call the  {\it discrepancy\/}. Namely, we define
   \[
\Delta_A(f,N) =     \max_{\psi\in \varPsi_N(A)}  
\left|\left \langle \OPN(f)\psi, \psi\right \rangle -  \int_{\TT^{2d}}
  f(\vx)d\vx  \right| ,
\] 
where    $ \varPsi_N(A)$ is the set of all normalised (that is, with $\|\psi\|=1$)
 eigenfunctions $\psi$ of $U_N(A)$ in $\HN$.

 Then the equidistribution property  for $N$ running through a certain infinite sequence 
 $\cN \subseteq \N$ means that  for all $f\in C^\infty(\TT^{2d})$, 
   \begin{equation}\label{QUE for all N's}
\lim_{\substack{N\to \infty \\ N \in \cN }}  \Delta_A(f,N)   = 0 . 
\end{equation}

\subsection{Previous results and methods}
The problem of bounding the discrepancy $\Delta_A(f,N)$  originates from the work of 
Kurlberg and Rudnick~\cite{KR2001},  establishing~\eqref{QUE for all N's} for $d=1$ 
and for  a set  $\cN$ of asymptotic density $1$. The method of~\cite{KR2001} combines 
some classical approaches to bounding {\it Kloosterman sums}, used before the {\it Weil 
bound\/} (see~\cite[Theorem~11.11]{IwKow}) became available, and some delicate results from 
the so-called {\it anatomy of integers\/}.

Bourgain~\cite{Bou2}, using bounds of exponential sums which 
rely on methods of additive combinatorics, 
 has improved the result of~\cite{KR2001} as 
  \begin{equation}\label{eq: Bourg}
    \Delta_A(f,N) \le N^{-\delta},
 \end{equation}
  for some unspecified $\delta > 0$ and also  for almost all $N$. 

Finally,  still for $d=1$, Ostafe, Shparlinski and Voloch~\cite{OSV}, using some classical methods of algebraic geometry, have shown that  one can take any $\delta< 1/60$ in~\eqref{eq: Bourg}. 

For $d\ge 2$, the only known 
result is due to Kurlberg, Ostafe,   Rudnick and Shparlinski~\cite{KORS}, which uses a combination 
of ideas of~\cite{Bou2} and~\cite{KR2001},   and gives~\eqref{QUE for all N's} in any dimension 
for a set $\cN$ of density $1$.

Note that while no rate of convergence is given in the works~\cite{KORS,KR2001}, 
this can be extracted from their arguments. However, these arguments are incapable of 
giving a power saving as in~\eqref{eq: Bourg}. 

Our previous work~\cite{BhaShp} introduces another point of view: instead of 
establishing~\eqref{QUE for all N's}, or even~\eqref{eq: Bourg} for almost all integers $N$, 
one can ask for explicit sequences of $N$ (preferably not growing too fast) 
for which~\eqref{eq: Bourg} holds for any $d \ge 1$ 
and also with a reasonably large value of $\delta$. 
Towards this question, a construction of~\cite{BhaShp} shows that for any $d$ and $A$
(as in the above), for 
a fixed but specially selected prime $p$, the sequence 
\begin{equation}
\label{eq:Seq pk}
\cN=\{p^k:~k =0,1, \ldots\}.
\end{equation}
satisfies 
\begin{equation}
\label{eq:B-S bound}
    \Delta_A(f,N) \le N^{-\kappa_d + o(1)},
\end{equation}
where 
\[
 \kappa_1 = 1/4, \qquad  \kappa_2 = 1/7,
\] 
 and for $d \ge 3$ we set
 \begin{equation}
\label{eq:kappa_d}
 \kappa_d =  \begin{cases} 
   \dfrac{3d}{2\(6d^2-5d+6\)},& \text{if} \ d \equiv 0 \pmod 3,\\
  \dfrac{3d^2 -1 }{2d\(6d^2 - 5d + 5\)},& \text{if} \ d \equiv 1 \pmod 3,\\
 \dfrac{3d}{2\(6d^2 - 5d +7\)},& \text{if} \  d \equiv 2 \pmod 3.
\end{cases}
\end{equation}
The approach of~\cite{BhaShp}  is based on a bound from~\cite{Shp}  
on exponential sums with linear recurrence sequences modulo large powers of a fixed prime.

We also mention the work of Kelmer~\cite{KelmerAHP} on prime power moduli $N =p^k$.
However, it applies to a slightly different scenario and also assumes that $k$ is fixed, but 
the prime $p$ grows.

\subsection{New discrepancy bound} 
Here we modify the argument of~\cite{BhaShp} and the way the bound of~\cite{Shp}  
is used,  and improve~\eqref {eq:B-S bound} by essentially doubling the decay rate, 
see~\eqref{eq:lim}. This new argument 
stems from the inductive idea of Bourgain~\cite{Bou1}, which is based 
on the interplay between exponential sums and counting solutions to the systems of congruences. 
It exploits the observation that discarding some congruences increases the main term, but gains in the quality of bounds on exponential sums, which control the error term. Proceeding this way recursively, with an optimal trade-off at each step, leads to our bound.

 We first recall the construction from~\cite{BhaShp}, 
 where we always assume that $A \equiv I_{2d}\pmod 2$ and that the characteristic polynomial  $F_A$ of the matrix $A\in \Sp(2d,\Z)$ is irreducible over $\Z$. 
This implies that  $A$ is diagonalisable over $\C$. 

We additionally 
assume that there are {\it no roots of unity\/} amongst the eigenvalues 
$\lambda_1, \ldots, \lambda_d, \lambda_{d+1} = \lambda_1^{-1},
 \ldots \lambda_{2d} = \lambda_d^{-1}$ of $A$ and their 
non-trivial ratios $\lambda_i/\lambda_j $, $1 \le i < j \le 2d$. 

By the Chebotarev density theorem, the primes $p > 2d$
for which  $F_A$ 
splits completely into distinct linear factors modulo $p$ form a set of 
 positive relative asymptotic density in the set of all primes. 
 Fix any such prime $p$, which clearly can be effectively found
 in a finite number of steps via a trial-and-error algorithm.  

Then it is shown in~\cite{BhaShp} that for this choice of $p$ 
and the sequence~\eqref{eq:Seq pk}, 
we have the bound~\eqref{eq:B-S bound}.

Our main result establishes~\eqref{QUE for all N's} for the above  sequence $\cN$ with 
a reasonably strong bound on rate of convergence. 

To formulate our result, we  set $\gamma_{1,s} = 0$ and for integers $e\ge 2$ and $h\ge 1$ we define 
 \begin{equation}
\label{eq:gamma_eh}
\gamma_{e,h} = 
 \min_{\substack{g_1,\ldots,g_{e-1}\in\Z\\
 g_1,\ldots,g_{e-1}\ge 0\\
                    g_1+\cdots+g_{e-1} =h-1}}\max\left\{0, 
\max_{1\leq \nu \leq e-1}
                   \sum_{j=\nu}^{e-1}\(1 -2 g_j/(j+1) \)\right\}, 
\end{equation}
 where the outer minimum is taken over all possible choices
 of non-negative integers $g_1,\ldots,g_{e-1}$ with $g_1+\cdots+g_{e-1} =h-1$.
 
 Next, we define 
 \begin{equation}
\label{eq:eta_d}
\eta_d = \max_{s =1, 2, \ldots } \frac{d -  \gamma_{2d,s}}{2s} .
\end{equation}  

\begin{theorem}\label{thm:DiscrBound}
For the above choice of prime $p$, for the sequence $\cN$ given by~\eqref{eq:Seq pk} and $N \in \cN$, 
we have 
\[ 
\Delta_A(f,N)  \le  N^{-\eta_d+o(1)}, 
\] 
as $N\to \infty$, where $\eta_d$ is given by~\eqref{eq:eta_d}.
  \end{theorem}

  \subsection{Numerical and asymptotic behaviour of the new rate of decay}
  \label{sec:num asymp}
We compare the previous exponent $\kappa_d$ in~\eqref{eq:kappa_d} 
with the new exponent $\eta_d$ in~\eqref{eq:eta_d}. 

First, we consider their numerical values. Table~\ref{tab:kappa-vs-eta} compares 
 $\kappa_d$  and $\eta_d$ for $d=2, \ldots, 10$;
 for $d=1$ we have $\kappa_1=\eta_1 = 1/4$. 

\begin{table}[htbp]
\centering
\setlength{\tabcolsep}{4pt}
\renewcommand{\arraystretch}{1.8}
\begin{tabular}{l*{10}{c}}
\toprule
$d$ &  $2$ & $3$ & $4$ & $5$ & $6$ & $7$ & $8$ & $9$ & $10$ \\
\midrule
$\kappa_d$ &   $\frac{1}{7}$ & $\frac{1}{10}$ & $\frac{47}{648}$ & $\frac{5}{88}$ & $\frac{3}{64}$ & $\frac{73}{1848}$ & $\frac{4}{117}$ & $\frac{9}{298}$ & $\frac{299}{11100}$ \\
$\eta_d$ &   $\frac{1}{6}$ & $\frac{43}{330}$ & $\frac{2}{19}$ & $\frac{11}{126}$ & $\frac{37}{495}$ & $\frac{4442}{68211}$ & $\frac{25642}{444015}$ & $\frac{243697}{4703985}$ & $\frac{5177218}{110160765}$ \\
$\eta_d/\kappa_d$ &  $1.17$ & $1.30$ & $1.45$ & $1.54$ & $1.59$ & $1.65$ & $1.69$ & $1.72$ & $1.74$ \\
\bottomrule
\end{tabular}
\caption{Values of $\kappa_d$ and $\eta_d$ and approximate ratios.}
\label{tab:kappa-vs-eta}
\end{table}

We now show that for every integer $d\ge1$,
\begin{equation}\label{eq:etad large}
\eta_d >
 \frac{d}{2d^2+d+2H_{2d}}
\end{equation}
where 
\[
H_n=\sum_{k=1}^n k^{-1}
\]

We   give an explicit the choice of $g_1, \ldots, g_{2d-1}$ 
which leads to the desired lower bound~\eqref{eq:etad large}.

Set  $\sigma_{2d}=0$ and then for $j=2d-1,\ldots,1$ define recursively (in the descending 
order)
\begin{equation}\label{eq:recurs}
g_j=\rf{\frac{(j+1)\(1+\sigma_{j+1}\)}{2}} 
\quad \text{and} \quad 
\sigma_j=\sigma_{j+1}+1-\frac{2g_j}{j+1}.
\end{equation}
In particular, 
\begin{equation}\label{eq:residual}
0  \ge \sigma_j> -\frac{2}{j+1},  
\end{equation} 
and thus $g_j\ge0$ for every $j =1, \ldots, 2d-1$. 

We see from~\eqref{eq:recurs} and~\eqref{eq:residual} that,  
for $\nu =1, \ldots, 2d-1$, we have 
\begin{equation}\label{eq:negative}
\sum_{j = \nu}^{2d-1}  \(1-\frac{2g_j}{j+1} \) 
= \sum_{j = \nu}^{2d-1}  \( \sigma_j - \sigma_{j+1} \)
= \sigma_\nu - \sigma_{2d}  = \sigma_\nu  \le 0.
\end{equation} 
Therefore, choosing 
\[
s = g_1+\ldots+g_{2d-1}  + 1
\]
we conclude from~\eqref{eq:negative} that $\gamma_{2d,s} =0$. Since 
\[
g_j=\frac{j+1}{2}\(1+\sigma_{j+1}-\sigma_j\), \qquad  j =1, \ldots, 2d-1,
\]
using summation by parts, we derive
\[
s-1=\frac{(2d-1)(2d+2)}{4}
       -\sigma_1-\frac{1}{2}\sum_{j=2}^{2d-1}\sigma_j.
\]

Recalling~\eqref{eq:residual}, we obtain 
\[
s - 1 <\frac{(2d-1)(d+1)}{2} 
  + 1 + \sum_{j=2}^{2d-1}\frac{1}{j+1}
=\frac{2d^2+ d }{2}+H_{2d} - 1, 
\]
which now implies~\eqref{eq:etad large}. 

While it is not the main point, we now show that~\eqref{eq:etad large} is close to an optimal bound.
Indeed, for each $s =1,2, \ldots$,  choose an optimal tuple 
$\(g_1, \ldots g_{2d-1}\)$,  attaining $\gamma_{2d,s}$.  After this, we  set 
 $\sigma_{2d}= 0$ and 
\[
\sigma_\nu = \sum_{j=\nu}^{2d-1}\(1 -2 g_j/(j+1) \), \qquad \nu =1, \ldots, 2d-1.
\]
Thus  $g_j = (j+1)\(1+ \sigma_{j+1}- \sigma_j \)/2$, $j=1, \ldots, 2d-1$
\begin{align*}
s-1 &  = \sum_{j=1}^{2d-1} g_j  = \frac{d(2d+1)-1}{2}+ \frac{1}{2} 
\sum_{j=1}^{2d-1} (j+1)\(\sigma_{j+1}- \sigma_j \)\\
& =   \frac{2d^2+d-1}{2}+ d \sigma_{2d} -  \sigma_1 - \frac{1}{2} \sum_{j=2}^{2d-1} \sigma_j.
\end{align*}
From the trivial inequality $\sigma_\nu \le \gamma_{2d,s}$, we now conclude that 
\[
 s \ge \frac{2d^2+d+1}{2} - d \gamma_{2d,s}.
 \]
 Using this inequality with $s$ that corresponds to $\eta_d$ we derive 
\begin{equation}\label{eq:etad not too large}
 \eta_d \le   \frac{d} {2d^2+d+1} .
 \end{equation}
 In particular, we see from~\eqref{eq:kappa_d}, \eqref{eq:etad large} and~\eqref{eq:etad not too large} that 
\begin{equation}\label{eq:lim}
\lim_{d \to \infty} \eta_d/\kappa_d = 2.
 \end{equation}
 
 \section{Preliminaries} 

 \subsection{Notation}  
 Throughout the paper, the  notations 
\[X = O(Y),  \qquad X \ll Y, \qquad Y \gg X
\] 
are all equivalent to the
statement that the inequality $|X| \le c Y$ holds with some 
constant $c> 0$, which may depend on the matrix $A$.

We also write 
\[
X \asymp Y
\]
to denote $X \ll Y \ll X$.

 We recall that the additive character with period $1$ is denoted by
\[
z \in \mathbb R\  \mapsto \ \e(z) =  \exp\(2\pi i z\).
\]
For an integer $q \ge 1$ it is also convenient to define 
\[
\eq(z) = \e(z/q).
\]

The letter $p$, with or without indices, always denotes prime numbers.

Given an algebraic number $\gamma$ we denote by $\ord(\gamma, N)$ its order modulo $N$ 
(assuming that the ideals generated by $\gamma$ and $N$ are relatively prime in an appropriate number field). In particular, for an element $\lambda \in\F_{p^s}$, $\ord(\lambda,p)$ represents the order of $\lambda$ in $\F_{p^s}$. 

Similarly, we use $\ord(A,N)$ to denote the order of $A$ modulo $N$ (which always exists 
if $\gcd(\det A, N)= 1$ and in particular for $A\in  \Sp(2d,\Z)$). 

Finally, we use $\nu_p(z)$ to denote the $p$-adic order of $z \in \Q_p$,  where $\Q_p$
is the field of  $p$-adic numbers.

\subsection{Discrepancy and operators  $\TN$} 
 \label{sec:Discr TN}
  As in~\cite{KR2001}, and then also in~\cite{Bou2,KORS,OSV}, we observe that 
it is enough to bound  the quantity 
$\langle \TN(\vu) \psi, \psi \rangle$,
where
\[
\TN(\vu)=  \OPN(\e\( \vx\cdot \vu \)), 
\]
for any nonzero vector  $\vu \in\Z^{2d}$, 
see also~\eqref{eq:OpN}, and $\psi \in \Psi_N(A)$ runs through  eigenfunctions 
of $U_{N}(A)$, with frequency $\vu$ growing sufficiently slowly with $N$ (for example, as 
a power $N^\tau$ for a sufficiently small $\tau> 0$).

In fact, taking into account the rapid decay of the Fourier coefficients of functions  
$f\in C^\infty(\TT^{2d})$, we  see that to achieve the bound of 
Theorem~\ref{thm:DiscrBound} it is enough to prove that for $N \in \cN$, 
for all non-zero vectors $\vu \in \Z^{2d}$, we have
\begin{equation}
\label{eq:TN bound}
   \max_{\psi \in \Psi_N(A)} 
  \left| \langle \TN(\vu)\psi ,\psi  \rangle \right| \le C \|\vu\|^C N^{-\eta_d + o(1)},
\end{equation}
where $\|\vu\|$ is the Euclidean length of $\vu$,  
with some constant $C> 0$ depending only on $A$ and $p$ and with $\eta_d$ given by~\eqref{eq:eta_d}.

\subsection{Operators  $\TN$ and  congruences} 
 \label{sec:TN congr}
 In turn, the proof of~\eqref{eq:TN bound} is based on a link  between $ \left| \langle \TN(\vu)\psi ,\psi  \rangle \right| $ and bounds on the 
number of solutions to certain systems of congruences, first established in~\cite{KR2001}
and then generalised and used in all other papers on this 
subject~\cite{Bou2,KORS,OSV}. 

By our choice of $p$, the eigenvalues of $A$ are pairwise distinct modulo $p$.  Hence, 
we can assume,  using Hensel lifting, that the eigenvalues $ \lambda_{1},\ldots,\lambda_{2d}$ of $A$ 
are elements of the reduced residue system $ \Z_{p^k}^*$ modulo $p^k$, that is, as in~\cite{BhaShp} we can assume  that 
\[ 
 \lambda_{1},\ldots,\lambda_{2d}\in  \Z_{p^k}^*. 
\]

Namely,  let  $ R_{r,s}$ be the number of solutions to the system of  congruences
\[
 \lambda_{i}^{x_1}+\ldots + \lambda_{i}^{x_{s}}\equiv \lambda_{i}^{y_1}+\ldots + \lambda_{i}^{y_{s}} \pmod {p^{r}}, \quad i=1,\ldots,2d, 
\]
in variables $x_{1}, y_{1}, \ldots, x_{s}, y_{s} =1, \ldots, \ord(A,p^r)$.  
Then, combining~\cite[Lemma~2.1]{BhaShp} and~\cite[Equation~(2.5)]{BhaShp}, 
we see that for  any fixed $s$, if $k > m$ then  for any nonzero vector  $\vu \in\Z^{2d}$  we have 
\begin{equation}
\label{eq:TN and R}
   \max_{\psi \in \Psi_N(A)} 
  \left| \langle \TN(\vu)\psi ,\psi  \rangle \right|  ^{2s}
   \ll p^{dk+ 2ms} \frac{R_{k-m,s}}{\ord(A,p^k)^{2s} }, 
\end{equation}
where $N = p^k$, $m\ge 1$ is defined by 
\begin{equation}
\label{eq:pm > u}
p^{m} \ge c \|u\| ^{2d} > p^{m-1}
\end{equation}
for some constant $c > 1$, depending only on $A$ (as in~\cite[Lemmas~2.3 and 2.4]{BhaShp}).

  \subsection{Multiplicative orders}  
We recall some results from~\cite{BhaShp}, which in turn are variations 
of results of Korobov~\cite{Kor1,Kor2}.

First we recall the following well-known property of orders; see, 
for example,~\cite[Lemma~3.1]{BhaShp}.

  \begin{lem}\label{eq:ord lin growth} 
Assume that a prime  $p\ge 3$ and an integer  $\lambda\ne \pm 1$ are relatively prime. 
Let 
\[
\gamma = \nu_p\(\lambda^{\ord(\lambda, p)} -1\).
\]
Then for any integer  $r \ge \gamma$ we have 
\[
\ord(\lambda, p^r) = \ord(\lambda, p) p^{r-\gamma}.
\]
\end{lem}

In particular, combining Lemma~\ref{eq:ord lin growth}  with~\cite[Lemma~3.2]{BhaShp},  
we derive: 
 
 \begin{lem}\label{eq:ord growth} 
There is a constant $\xi> 0$ depending only on $p$ and $A$, such that 
\[
p^r\ge \ord(\lambda_i, p^r), \ord(\lambda_i/\lambda_j, p^r) \ge \xi p^{r}, 
\qquad 1 \le i, j \le 2d, \ i \ne j.
\]
\end{lem}

 \subsection{Exponential sums}
The following bound on these exponential sums follows from more general results 
of~\cite{Shp}, it  is also stated in~\cite[Lemma~3.3]{BhaShp} with a sketch of a proof. 

\begin{lem}\label{lem:exp sum}
Fix an integer $e \ge 1$ and let $\mu_1, \ldots,  \mu_{e}\in \Z_{p^r}^*$ be pairwise distinct modulo $p$ elements  
and such that 
\[
\ord(\mu_i, p^r), \ord(\mu_i/\mu_j, p^r) \ge \xi p^{r}, 
\qquad 1 \le i, j \le e, \ i \ne j, 
\]
for some constant $\xi > 0$. 
  Then,  for 
 \[
S_r(\va, \bmu) =   \sum_{x=1}^{t_r(\va, \bmu)}  \e_{p^r}\(a_1\mu_1^x + \ldots + a_{e}\mu_{e}^x \), 
\]
where  $\va= (a_1, \ldots, a_e)$, $\bmu = \mu_1, \ldots,  \mu_{e}$ and 
 $t_r(\va, \bmu)$   is the period of the sequence $
 a_1\mu_1^x + \ldots + a_{e}\mu_{e}^x$, $x =0,1, \ldots$
 uniformly over $\bmu$   and  $\va$  with 
$\gcd(\va, p) =  1$, 
we have
\[
\left| S_r(\va, \bmu)  \right| \le 
t_r(\va, \mu)^{1-1/e+o(1)}, \qquad \text{as} \ r \to \infty.
 \]
\end{lem}

\begin{rem} In our application of Lemma~\ref{lem:exp sum} we always have 
\[
\{\mu_1, \ldots,  \mu_{e}\} \subseteq \{\lambda_1, \ldots,  \lambda_{2d}\},
\] 
however, examining the argument in~\cite{BhaShp, Shp}, one can verify that it holds with 
the claimed uniformity with respect to $\bmu$ (provided $\bmu$ satisfies the improsed conditions 
on $\ord(\mu_i, p^r)$ and $\ord(\mu_i/\mu_j, p^r)$. 
\end{rem}

 \section{Proof of Theorem~\ref{thm:DiscrBound}}

\subsection{Bounding  the number of solutions to systems of exponential equations}   
To estimate $ R_{r,s}$ from Section~\ref{sec:TN congr}, 
we consider a more general question. 

Given $\bmu = \(\mu_1, \ldots,  \mu_{e}\)$ with $e \ge 1$   pairwise distinct modulo $p$ elements  
$\mu_1, \ldots, \mu_{e}\in \Z_{p^r}^*$, we denote 
by  $  Q_{r,s}\(\bmu \)$ the number of solutions  to the system of  congruences
\[
 \mu_{i}^{x_1}+\ldots + \mu_{i}^{x_{s}}\equiv \mu_{i}^{y_1}+\ldots + \mu_{i}^{y_{s}} \pmod {p^{r}}, \quad i=1,\ldots,e, 
\]
in variables $x_{1}, y_{1}, \ldots, x_{s}, y_{s} =1, \ldots,  \ell_r(\bmu)$, 
where 
\[
 \ell_r(\bmu) = \lcm\left[\ord\(\mu_1,p^r\), \ldots, \ord\(\mu_e,p^r\)\right].
 \]
 
 It is useful to observe that for $\va= (a_1, \ldots, a_e)$ 
 with $\gcd(\va, p) =  1$
 the period   $t_r(\va, \bmu)$, defined in Lemma~\ref{lem:exp sum} 
satisfies 
\[
 \min\left\{\ord\(\mu_1,p^r\), \ldots, \ord\(\mu_e,p^r\)\right\} \le t_r(\va, \bmu) \mid   \ell_r(\bmu) 
\]
Indeed, the divisibility is trivial. To see the lower bound, we observe that 
\[
z_1 = a_1\(\mu_1^{t_r(\va, \bmu)} -1\), \ldots,  z_e = a_{e} \(\mu_e^{t_r(\va, \bmu)} -1\) 
\]
is a solution to 
 the following system of  $e$ homogeneous congruences modulo $p^r$  in $e$ variables
\[
\mu_1^i z_1 + \ldots +  \mu_{e}^i z_e
\equiv 0 \pmod{p^r}, \qquad 
 i = 0, \ldots, e-1.
\]

The Vandermonde
determinant of this system does not vanish modulo $p$ since by our assumption 
that $\mu_1, \ldots,  \mu_{e}$  are  pairwise distinct modulo $p$. 
Hence, $z_j \equiv 0 \pmod{p^r}$ for  all  $j =1, \ldots, e$.
Since  $\gcd(\va, p) = 1$, for at least one $j =1, \ldots, e$, we have 
$\mu_j^{t_r(\va, \bmu)} -1\equiv 0 \pmod {p^r}$. Therefore, 
we have $\ord\(\mu_j,p^r\) \mid t_r(\va, \bmu)$ 
and thus the desired lower bound on $t_r(\va, \bmu)$. 

In particular (as in Lemma~\ref{lem:exp sum}), throughout our argument, we use  that if 
\[
\ord(\mu_i, p^r), \ord(\mu_i/\mu_j, p^r) \ge \xi p^{r}, 
\qquad 1 \le i, j \le e, \ i \ne j, 
\]
then 
for any $\va= (a_1, \ldots, a_e)$  with $\gcd(\va,p)= \gcd(a_1, \ldots, a_e, p) =  1$, we have 
  \begin{equation}
 \label{eq: t and ell}
t_r(\va, \bmu) \asymp  \ell_r(\bmu) \asymp p^r.
\end{equation}
with the implied constants depending only on $p$ and $\xi$.  

Since the eigenvalues $\lambda_1, \ldots, \lambda_{2d}$ of $A$ are pairwise distinct 
modulo $p$,  it is easy to show that 
\[
\ord(A,p^r) = \ell_r(\blambda),
\]
where $\blambda = \(\lambda_1, \ldots, \lambda_{2d}\)$. 
Therefore,  
  \begin{equation}
 \label{eq: Q and R}
 Q_{r,s}\(\blambda\) =  R_{r,s}
\end{equation}
 where $ R_{r,s}$ is as in Section~\ref{sec:TN congr}.
 
 It is convenient to set $Q_{0,h}\(\bmu\) = 1$.

\begin{lem}\label{eq:recurs Q} 
Let $\mu_1, \ldots,  \mu_{e}\in \Z_{p^r}^*$  be pairwise distinct modulo $p$ and also such that 
\[
\ord(\mu_i, p^r), \ord(\mu_i/\mu_j, p^r) \ge \xi p^{r},
\qquad 1 \le i, j \le e, \ i \ne j, 
\]
for some constant $\xi > 0$. Then,  for any  
integers $g\ge 0$ and $h\ge 1$ with $g+h = s$ we have 
 \[
    Q_{r,s}\(\bmu\) 
\le    p^{r(2s-e+o(1))} \sum_{u = 0}^{r} p^{eu -2hu - 2gu/e} 
 Q_{u,h}\(\bmu\). 
\]
\end{lem} 

 \begin{proof} By orthogonality of additive characters, we have
  \begin{equation}
 \label{eq:R = S(a)}
    Q_{r,s}\(\bmu \) 
   =  \frac{1}{p^{er}}\sum_{\va \in \Z_{p^r}^{e}}\left|\widetilde S_r(\va, \bmu) \right|^{2s},
\end{equation}
where 
  \begin{equation}
 \label{eq:S-bar S}
\widetilde S_r(\va, \bmu) =    \sum_{x=1}^{\ell_r(\bmu)}  \e_{p^r}\(a_1\mu_1^x + \ldots + a_{e}\mu_{e}^x \)
=  \frac{\ell_r(\bmu)} {t_r(\va, \bmu)} S_r(\va, \bmu),
\end{equation}
where $S_r(\va, \bmu)$ is the exponential sum as defined in Lemma~\ref{lem:exp sum}.

Next, for each $u= 0, \ldots, r$,  we collect together $O\(p^{ue}\)$ vectors with 
\[
\gcd(\va, p^r) =  p^{r-u}
\]
and notice that for $u \ge 1$ we have
  \begin{equation}
 \label{eq: Sr Su}
\widetilde  S_r(\va, \bmu)  = \frac{\ell_r(\bmu)}{\ell_{u}(\bmu)} \widetilde S_{u}(\va p^{-r+u}, \bmu) .
\end{equation}
We also  see from Lemma~\ref{eq:ord lin growth} (or, alternatively, from~\eqref{eq: t and ell})  that 
\[
\frac{\ell_r(\bmu)}{\ell_{u}(\bmu)} \asymp p^{r-u}.
\]

Writing $\va = p^{r-u} \vb$ for $\vb = (b_1, \ldots, b_e) \in   \Z_{p^{u}}^{e}$ 
with $\gcd(\vb, p) = \gcd(b_1, \ldots, b_{e}, p) =  1$ and separating the contribution 
\[
  \frac{1}{p^{er}}   \ell_r(\bmu)^{2s}  \ll  p^{r(2s-e)}
\]
of $u=0$,  we  derive from~\eqref{eq:R = S(a)} and~\eqref{eq: Sr Su} that  
\begin{align*}
    Q_{r,s}\(\bmu\) & =\frac{1}{p^{er}} \ell_r(\bmu)^{2s}   +  \frac{1}{p^{er}} \sum_{u = 1}^{r} 
    \( \frac{\ell_r(\bmu)}{\ell_{u}(\bmu)}\)^{2s}    \sum_{\substack{\vb \in \Z_{p^{u}}^{e}\\ \gcd(\vb, p)=1}}\left| \widetilde   S_{u}(\vb, \bmu) \right|^{2s}\\
    &   \ll   p^{r(2s-e)} +  \frac{1}{p^{er}} \sum_{u = 1}^{r} p^{2s(r-u)} 
   \sum_{\substack{\vb \in \Z_{p^{u}}^{e}\\ \gcd(\vb, p)=1}}\left| \widetilde   S_{u}(\vb, \bmu) \right|^{2s}. 
\end{align*}
Next,  we write 
\[
\left|  \widetilde S_{u}(\vb, \bmu) \right|^{2s} = \left|  \widetilde S_{u}(\vb, \bmu) \right|^{2g}\left|  \widetilde S_{u}(\vb, \bmu) \right|^{2h}.
\]
After this, we apply the bound of Lemma~\ref{lem:exp sum} to the $2g$-th power of $\left|  \widetilde S_{u}(\vb, \bmu) \right| \asymp \left|   S_{u}(\vb, \bmu) \right|$, which is of exactly the same form by~\eqref{eq: t and ell} and~\eqref{eq:S-bar S}, 
and leave the rest unchanged. This  gives us 
\begin{align*}
   Q_{r,s}\(\bmu\) 
& \ll  p^{r(2s-e)} \\
& \qquad \quad +   \frac{1}{p^{er}} \sum_{u = 1}^{r}p^{2s(r-u)+2gu(1-1/e+o(1))} 
   \sum_{\substack{\vb \in \Z_{p^{u}}^{e}\\ \gcd(\vb, p)=1}}\left| \widetilde S_{u}(\vb, \bmu) \right|^{2h}\\
& = p^{r(2s-e)} +  p^{r(2s-e+o(1))} \sum_{u = 1}^{r}p^{-2hu - 2gu/e} 
   \sum_{\substack{\vb \in \Z_{p^{u}}^{e}\\ \gcd(\vb, p)=1}}\left| \widetilde S_{u}(\vb, \bmu) \right|^{2h}.
\end{align*}

Finally, extending the summation over $\vb$ with $\gcd(\vb, p)=1$ to all vectors $\vb \in \Z_{p^{u}}^{e}$ and 
recalling the orthogonality relation~\eqref{eq:R = S(a)}, we obtain 
\[
    Q_{r,s}\(\bmu\) 
\le  p^{r(2s-e)} +   p^{r(2s-e+o(1))} \sum_{u = 1}^{r} p^{eu -2hu - 2gu/e} 
 Q_{u,h}\(\bmu\). 
\]
Recalling our convention that $Q_{0,h}\(\bmu\) = 1$, we conclude the proof. 
\end{proof}

\begin{lem}\label{lem: bound Q} 
Let $\mu_1, \ldots,  \mu_{e}\in \Z_{p^r}^*$ be pairwise distinct modulo $p$ 
 such that 
\[
\ord(\mu_i, p^r), \ord(\mu_i/\mu_j, p^r) \ge \xi p^{r},
\qquad 1 \le i, j \le e, \ i \ne j, 
\]
for some constant $\xi > 0$. 
  Then,  
 \[
    Q_{r,s}\(\bmu\) 
\le    p^{r(2s-e+\gamma_{e,s} + o(1))},
\]
where $\gamma_{e,s}$ is given by~\eqref{eq:gamma_eh}.
\end{lem}

 \begin{proof} We prove the desired estimate by induction on $e$. 
 For $e = 1$, that is, $\bmu = \mu \in \Z_{p^r}^*$ 
 we trivially have $Q_{r,s}\(\mu\) \le p^{r(2s-1)}$, 
 while $\gamma_{1,s} = 0$. 
 
 Now assume the desired bound holds for $e-1$ (any  modulus $p^u$, $u =1, 2, \ldots$). 
 Denote $\bmu^* = \(\mu_1, \ldots, \mu_{e-1}\)$. 

Using the trivial inequality 
 \[
     Q_{r,s}\(\bmu\)  \le  \(\frac{\ell_r\(\bmu\)}{\ell_r\(\bmu^*\)}\)^{2s}  Q_{r,s}\(\bmu^*\) 
     \ll Q_{r,s}\(\bmu^*\).
 \]
 we derive from Lemma~\ref{eq:recurs Q} that for any
integers $g\ge 0$ and $h\ge 1$ with $g+h = s$
\[
   Q_{r,s}\(\bmu\)  \le p^{r(2s-e+o(1))} \sum_{u = 0}^{r} p^{u(e -2h - 2g/e)}  Q_{u,h}\(\bmu^*\).
\]
Applying the inductive assumption, we now see that 
\begin{align*}
   Q_{r,s}\(\bmu\) & \le   p^{r(2s-e+o(1))} \sum_{u = 0}^{r} p^{u(e -2h - 2g/e)}  
\(1+ p^{u(2h-e +1+\gamma_{e-1,h})}\)\\
& \le p^{r(2s-e+o(1))} \(1 + p^{r(e -2h - 2g/e)}   + p^{r(\gamma_{e-1,h} +1- 2g/e)}\).
\end{align*}
 Observe that for $e > 2$, for any admissible tuple $\(g_1, \ldots g_{e-2}\)$ 
in~\eqref{eq:gamma_eh} (with $e-1$ instead of  $e$) we have
\[
\sum_{j=1}^{e-2}\(1 -2 g_j/(j+1)\) \ge e-2 - \sum_{j=1}^{e-2}g_j = e-h-1
\]
and thus 
  \begin{equation}
 \label{eq:  gamma vs e-h-1}
 \gamma_{e-1,h}  \ge e-h -1 .
\end{equation}
For $e=2$ the above inequality also holds as $ \gamma_{e-1,h}  = 0 \ge 1 - h$.
 
Hence 
\[
 Q_{r,s}\(\bmu\)  \le p^{r(2s-e+o(1))} \(1 +   p^{r(\gamma_{e-1,h} +1- 2g/e)}\) .
 \]
We now choose $g$ and $h$  to realise the minimum in 
\[
 \gamma_{e,s}  = \min_{\substack{g \ge 0, \ h \ge 1\\ g+h = s}} \max\{0, \gamma_{e-1,h} +1- 2g/e\} ,
\]
which concludes the proof. 
\end{proof}
 
\subsection{Concluding the proof}  Our goal is to establish~\eqref{eq:TN bound}.
Since we have the trivial bound 
\[
   \max_{\psi \in \Psi_N(A)} 
  \left| \langle \TN(\vu)\psi ,\psi  \rangle \right| \le 1
 \]
 on the unitary operator $\TN(\vu)$, taking the constant $C$ in~\eqref{eq:TN bound}
 to be sufficiently large (in fact $C > 2d \eta_d$ suffices), we see that it holds for $k \le m$, where
 $m$ is defined by~\eqref{eq:pm > u}.
  
 So, we can now assume that $k> m$.

We recall the relation~\eqref{eq: Q and R} and observe that  by Lemma~\ref{eq:ord growth}, the bound of 
Lemma~\ref{lem: bound Q}  applies to $R_{k-m,s}$ 
and yields the bound
\[
R_{k-m,s} \le p^{(k-m)(2s-2d+ \gamma_{2d,s}) + o(k)} \le  p^{k(2s-2d+ \gamma_{2d,s}  + o(1))}
\]
(note that the second inequality holds since by~\eqref{eq:  gamma vs e-h-1} we have
$\gamma_{2d,s}  \ge 2d-s$). 

We now substitute this estimate into~\eqref{eq:TN and R}, and using that $p^m \ll \|\vu\|^{2d}$,  derive
\[
   \max_{\psi \in \Psi_N(A)} 
  \left| \langle \TN(\vu)\psi ,\psi  \rangle \right|  ^{2s} 
   \ll \|\vu\|^{C_0} p^{dk} \frac{p^{k(2s-2d+ \gamma_{2d,s}+ o(1))}}{\ord(A,p^k)^{2s} }, 
\]
where $C_{0}> 0$ depends only on $A$, $p$ and $s$ and can easily be evaluated 
explicitly (which is not necessary for our argument). 
Since by Lemma~\ref{eq:ord growth},  $\ord(A,p^k) \gg p^k$, we now conclude that 
\[
   \max_{\psi \in \Psi_N(A)} 
  \left| \langle \TN(\vu)\psi ,\psi  \rangle \right|  ^{2s} 
   \ll  \|\vu\|^{C_0} p^{-(d -  \gamma_{2d,s}+o(1))k} .
\]
Optimising $s$ leads to the bound~\eqref{eq:TN bound} and concludes the proof.

\section*{Acknowledgement} 
During the preparation of this work, the authors were   supported in part by the Australian Research Council Grant DP230100530. S.B. would like to thank IIIT-Delhi and IIT Bombay for their support and excellent working conditions during his visits while preparing this work. I.S. would also like to thank the Max Planck Institute for Mathematics, Bonn, for its hospitality, stimulating environment, and generous support. 

ChatGPT-6 has been used for language editing and for generation of numerical and asymptotic 
results in Section~\ref{sec:num asymp}.

\end{document}